\documentclass[10pt]{article}
\usepackage{latexsym,amssymb,amsmath,epsfig}
\usepackage{epsfig,graphicx,latexsym,graphics} 
\usepackage{epic}
\usepackage{amsfonts}
\usepackage{color}
\newtheorem{thm}{Theorem}
\newtheorem{lem}{Lemma}
\newtheorem{cor}{Corollary}

\newtheorem{con}{Conjecture}

\begin{document}


\centerline {\Large\bf{Extremal $\bf k$-ape Trees for Randi\'{c} Index }}

\centerline{ Naveed Akhter*, Muhammad Kamran Jamil, Ioan Tomescu}
\centerline {\footnotesize Abdus Salam School Of Mathematical Sciences}
\centerline {\footnotesize Government College University, Lahore Pakistan}
\centerline {\footnotesize akhtarnaweed@yahoo.com}

\begin{abstract}
The Randi\'{c} (connectivity) index is one of the most successful molecular descriptors in structure-property and structure-activity relationships studies. A graph $G$ is called an apex tree  if it contains a vertex $x$ such that $G-x$ is a tree. For any integer $k \ge 1$ the graph $G$ is called $k$-apex tree if there exists a subset $X$ of $V(G)$ of cardinality $k$ such that $G-X$ is a tree and for any $Y \subset V(G)$ and $|Y| < k$, $G-Y$ is not a tree. J. Gao \cite{g} found the sharp upper bound for the Randi\'{c} index of apex trees. In this paper, we proved that $k$-apex trees are not regular for $k\ge2$ and proposed a sharp upper bound for the Randi\'{c} index of $k$-apex trees for $k \ge 2$.
\end{abstract}

%

\section{Introduction}

In studying branching properties of alkanes, several numbering schemes for the edges of the associated hydrogen-suppressed graph were proposed based on the degrees of the end vertices of an edge \cite{MJ}. To preserve rankings of certain molecules, some inequalities involving the weights of edges needed to be satisfied. Randi\'{c} \cite{MJ} stated that weighting all edges uv of the associated graph $G$ by $(d(u)d(v))^{\frac{1}{2}}$ preserved these inequalities, where $d(u)$ and $d(v)$ are the degrees of vertices $u$ and $v$ respectively. The sum of weights over all edges of $G$, which is called the Randi\'{c}  index or molecular connectivity index or simply connectivity index of $G$ and denoted by $R(G)$, has been closely correlated with many chemical properties \cite{KH} and found to parallel the boiling point, Kovats constants, and a calculated surface. In addition, the Randi\'{c} index appears to predict the boiling points of alkanes more closely, and only it takes into account the bonding or adjacency degree among carbons in alkanes \cite{KH1}. It is said in \cite{HM} that Randi\'{c} index  together with its generalizations it is certainly the molecular-graph-based structure descriptor, that found the most numerous applications in organic chemistry, medicinal chemistry, and pharmacology. More data and additional references on the index can be found in \cite{GVN,GL}.

All graphs considered in this paper are simple, finite and connected. For a vertex $v \in V(G)$ its degree is denoted by $d_G(v)$ and if $G$ is clear from the context we  simplify the notation to $d(v)$. For $ X \subset V(G)$, $G-X$ is the subgraph of $G$ obtained from $G$ by removing the vertices of $X$ and edges incident with them, in particular $G-\{v\}$ is denoted by $G-v$. An edge $uv$ in $G$ is called \textit{symmetric} edge if $d(u)=d(v)$. An edge which is not symmetric is called \textit{asymmetric}.

The  Randi\'{c} index of the graph $G$ is defined as
\[R(G) = \sum\limits_{uv \in E\left( G \right)}{1 \over \sqrt {d\left( u \right)d\left( v \right)}} \]

In topological graph theory, graphs that contain a vertex whose removal yields a planar graphs play an important role, these graphs are called apex graphs \cite{A,M}. Along these lines a graph $G$ is called an apex tree \cite{KZK}  if it contains a vertex $x$ such that $G-x$ is a tree. The vertex $x$ is called apex vertex of $G$. Note that a tree is always an apex tree, hence a non-trivial apex tree is an apex tree which itself is not a tree. For any integer $k \ge 1$ the graph $G$ is called $k$-apex tree if there exists a subset $X$ of $V(G)$ of cardinality $k$ such that $G-X$ is a tree and for any $Y \subset V(G)$ and $|Y| < k$, $G-Y$ is not a tree. A vertex in $X$ is called $k$-apex vertex. Clearly, $1$-apex trees are precisely non-trivial apex trees. Apex trees and $k$-apex trees were introduced in \cite{KJH} under the name quasi-tree graphs and $k$-generalized quasi-tree graphs, respectively. J. Gao found the sharp upper bound of Randi\'{c} index for apex trees and in this paper, we give sharp upper bound on the Randi\'{c} index of $k$-apex trees for $k \ge 2$.

\subsection{Non-regularity of $k$-apex trees}

We need the following upper  bound on Randi\'{c} index to prove our main results.

\begin{lem}\cite{cghp}{\label{maxRan}}
    If $G$ is a connected graph of order $n$, then
    \[ R(G) \le {n \over 2} - {1 \over 2}{\sum\limits_{uv \in E( G )} \left( {1 \over {\sqrt {d( u )}} } - {1 \over {\sqrt {d( v )}}}  \right)^2} \]
\end{lem}

\begin{lem}{\label{incinx}}
    The following is an increasing function
    \[f(x)=\left({1 \over \sqrt{x}} - {1 \over \sqrt a} \right)^2 \] where $a$ and $x$ are positive real numbers and $x>a$.
\end{lem}

\begin{lem}{\label{decinx}}
    The function $f(x) = \left( {1 \over \sqrt{x+1} } - {1 \over \sqrt{x}} \right)^2 $ for $x > 0 $ is a decreasing function.
\end{lem}

\begin{thm}\label{notreg}
    If $G$ is a $k$-apex $(k \ge 2)$ tree of order $n \ge 4k-1$, then $G$ is not regular.
\end{thm}

\textbf{Proof.} Let $X \subset V(G)$ such that $ |X|=k$ and $G-X$ is a tree, then
\[ \sum_{v \in V(G-X)} d_{G-X}(v)=2n-2k-2 \]

Suppose that $G$ is an $m$-regular  graph. As $k$-apex  tree $(k \ge 2)$ is never two regular and pendant vertex in $G-X$ has at most $k+1$ degree in $G$, therefore   $ 3 \le m \le k+1$. Further suppose that $l$ is the number of edges in $G$ whose one end is in the set $V(G-X)$ and other is in the set $X$, then

\[ l \le mk \]

As $G$ is $m$-regular, therefore  $ \sum\limits_{v \in V(G-X)} d_{G}(v) = mn-mk $ and hence

\begin{equation*}
    \begin{split}
        l = & \sum\limits_{v \in V(G-X)} d_{G}(v) - \sum_{v \in V(G-X)} d_{G-X}(v) \\
          = & mn - mk - 2n + 2k+2
    \end{split}
\end{equation*}

As $l \le mk $, therefore  $ mn-mk-2n+2k-2 \le mk$ and hence
\[ m \le {2n-nk-2 \over {n-2k}} \]

As  $ 3 \le m $, therefore
 \[ 3 \le {2n-nk-2 \over {n-2k}} \]

 \[ n \le 4k-2 \]

 Which is contrary to our supposition that $n \ge 4k-1 $. Hence $G$ is not $m$-regular $( 2 \le m \le k+1)$ that is $G$ is not regular. $\hspace{1cm}\Box$

\section{Extremal $\bf {k}$-apex tress for Randi\'{c} index}

\begin{figure}[h!]
  \centerline{\includegraphics[width=10cm]{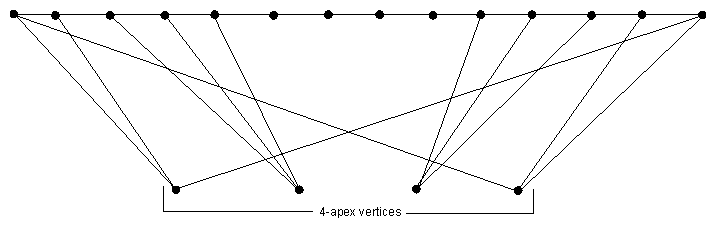}}
  \caption{$4$-apex tree of order 18.}\label{kapex4cycl}
\end{figure}

Suppose $\widetilde G_k^n$ is the set of all $k$-apex $(k \ge 2)$ trees of order $n$ $(n \ge 4k-1)$ such that each $H \in \widetilde G_k^n $ has vertices of degree two or three only and  has only two asymmetric edges. Then
\[ R(H)= \frac{n}{2}-\frac{5-2\sqrt 6}{6} \]
The Fig. \ref{kapex4cycl} is an example of a graph in $ \widetilde G_4^{18}$.

\begin{lem}
    If $a \ge 2$ is an integer and $ x \ge a+2 $, then the function
    \[ f(x) = \left( {1 \over \sqrt x}- \frac{1}{\sqrt a} \right)^2 - \left( \frac{1}{\sqrt 3}- \frac{1}{\sqrt2} \right)^2 \]
    is an increasing and positive function for $x\ge 4$.
\end{lem}

\begin{cor}\label{asymNTAlE}
    If $G$ is a $k$-apex tree of order $n \ge 4k-1$  which has asymmetric edges degrees of whose vertices are not almost equal, then for any $H \in \widetilde G_k^n$
    \[ R(G) < R(H) \]
\end{cor}

\begin{lem}
    If $x$ is a real number and $x \ge 4$ then
   \[f(x) = (x-1) \left( \frac{1}{\sqrt x} - \frac{1}{\sqrt{x-1}} \right)^2 - \left( \frac{1}{\sqrt 3}- \frac{1}{\sqrt 2} \right)^2 \]
    is greater than zero.
\end{lem}

\begin{lem}
    If $x$ is a real number and $x \ge m \ge 4$, then
    \[f(x) = (x-1)\left(\frac{1}{\sqrt m} - \frac{1}{\sqrt{m-1}} \right)^2 - \left( \frac{1}{\sqrt 3}- \frac{1}{\sqrt 2} \right)^2 \]
    is greater than zero.
\end{lem}

\begin{cor}\label{masymedg}
  If $G$ is a $k$-apex tree of order $n \ge 4k-1$, which has at least $2m-2$ $( 2 \le m \le k+2)$ asymmetric edges, degrees of whose vertices are almost equal then for any $H \in \widetilde {G}_k^n$
  \[ R(G) < R(H) \]
\end{cor}

\begin{figure}[h!]
  \centerline{\includegraphics[width=5cm]{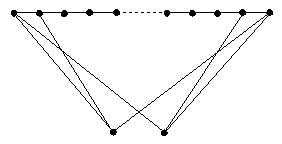}}
  \caption{A $2$-apex tree having sharp upper bound}\label{xx}
\end{figure}

\begin{con}\label{maxk-apextree}
    If $G$ is $k$-apex tree $(k \ge 2)$ of order $n \ge 4k-1$, then
    \[R(G) \le \frac{n}{2}- \frac{5-2\sqrt 6}{6} \]
    and the equality holds if and only if $ G \in \widetilde G_k^n$.
\end{con}


\section{Conclusion}

The Randi\'{c} index is one of the most successful molecular descriptors in structure-property and structure-activity relationships studies. J. Gao found the sharp upper bound for the Randi\'{c} index of apex trees. We found that for $k$-apex trees $(k \ge 2)$ of order $n \ge 4k-1$ no regular graphs exist and proposed a conjecture for sharp upper bound for Randi\'{c} index.

\end{document}